\documentclass[12pt]{article}
\usepackage[cp1251]{inputenc}
\usepackage[english]{babel}
\usepackage{mathrsfs}
\usepackage{amsfonts}
\usepackage[all]{xy}
\usepackage{amscd}
\usepackage{amsmath,amssymb,latexsym}
\usepackage{amsthm}

\newtheorem*{Pr1}{Proposition 1}
\newtheorem*{Pr2}{Proposition 2}
\newtheorem*{Pr3}{Proposition 3}
\newtheorem*{Pr4}{Proposition 4}
\newtheorem*{Pr5}{Proposition 5}
\newtheorem*{Pr6}{Proposition 6}
\newtheorem*{Pr7}{Proposition 7}
\newtheorem*{Pr8}{Proposition 8}
\newtheorem*{Pr9}{Proposition 9}
\newtheorem*{Le1}{Lemma 1}
\newtheorem*{Le2}{Lemma 2}

\newtheorem*{Cor}{Corollary 1}
\newtheorem*{Cor2}{Corollary 2}

\newtheorem*{th1}{Theorem 1}
\newtheorem*{TopPr}{Theorem 2}
\newtheorem*{def1}{Definition 1}
\newtheorem*{def2}{Definition 2}
\newtheorem*{thA}{Theorem A}
\newtheorem*{thB}{Theorem B}
\newtheorem*{thC}{Theorem C}

\def \id{{\rm id \,}}

\def \R{\mathbb{R}}
\def \N{\mathbb{N}}

\def\A{{\cal A}}
\def\S{{\cal S}}

\def \ss{\subset}
\begin{document}

\title{On Extension of Functors}
\author{L.Karchevska, T.Radul}
\date{\today}

\maketitle

Department of Mechanics and Mathematics, Lviv National University,
Universytetska st.,1, 79000 Lviv, Ukraine.\newline e-mail:
tarasradul@yahoo.co.uk, crazymaths@ukr.net



\begin{abstract} A.Chigogidze defined for each normal functor on the category $Comp$ an extension which is a normal functor on the category $Tych$. We consider this extension for any functor on the category $Comp$ and investigate which properties it preserves from the definition it preserves from the definition of normal functor. We investigate as well some topological properties of such extension. \end{abstract}

{\bf Introduction.} The general theory of functors acting in the
category $Comp$ of compact Hausdorff spaces (compacta) and
continuous mappings was founded by E.V.~Shchepin \cite{Shchepin}.
He distinguished some elementary properties of such functors and
defined the notion of normal functor that has become very
fruitful. The class of normal functors includes many classical
constructions: the hyperspace $\exp$, the functor of probability
measures $P$, the power functor and many other functors (see
\cite{ZT},\cite{FedZar} for more details). But some important
functors do not satisfy some of the properties from the Shchepin
list. Omitting some properties we obtain wider classes of functors
such as weakly normal functors and almost normal functors.

The properties from the definition of normal functor could be
easily generalized for the functors on the category  $Tych$ of
Tychonov spaces and continuous maps. Let us remark that $Tych$
contains $Comp$ as a subcategory.  A.Chigogidze defined for each
normal functor on the category $Comp$ an extension which is a
normal functor on the category $Tych$ \cite{ChExt}. This extension
could be considered for any functor on the category $Comp$. But
the situation is more complicated for wider classes of functors.
For example, the extension of the projective power functor (which
is weakly normal) does not preserve embeddings, which makes such
extension useless (see for example \cite{ZT}, p.67). However, if
we apply the Chigogidze extension to such  weakly normal functors
as the functor $O$ of order-preserving functionals, the functor
$G$ of inclusion hyperspaces, the superextension, we obtain
functors on the category  $Tych$ which preserve embeddings.

The main aim of this paper is to investigate which properties from
the definition of normal functor are preserved by Chigogidze
extension, specially we concentrate our attention on the
preserving of embeddings. The results devoted to this problem are
contained in Section 2. We define in this section the 1-preimages
preserving property which is crucial for preserving of embeddings.
In Section 3 we consider which functors have the 1-preimages
preserving property.

T.Banakh and R.Cauty obtained topological classification of the
Chigogidze extension of the functor of probability measures for
separable metric spaces. We generalize this result for convex
functors in Section 4.

{\bf 1.} All spaces are assumed to be Tychonov, all mappings are
continuous. All functors are assumed to be covariant. In the
present paper we will consider functors acting in two categories:
the category $Tych$ and its subcategory $Comp$.

Let us recall the definition of normal functor. A functor
$F:Comp\to Comp$ is called {\it monomorphic (epimorphic)} if it
preserves embeddings (surjections). For a monomorphic functor $F$
and an embedding $i:A\to X$ we shall identify the space $F(A)$ and
the subspace $F(i)(F(A))\ss F(X)$.

A monomorphic functor $F$ is said to be {\it preimage-preserving} if for each
map $f:X\to Y$ and each closed subset $A\ss Y$ we have $(F(f))^{-1}(F(A))=
F(f^{-1}(A))$.

For a monomorphic functor $F$ the {\it intersection-preserving} property is
defined as follows: $F(\cap\{X_\alpha\mid\alpha\in\A\})=\cap\{F(X_\alpha)\mid
\alpha\in\A\}$ for every family $\{X_\alpha\mid\alpha\in\A\}$ of closed subsets
of $X$.

A functor $F$ is called {\it continuous} if it preserves the limits of inverse
systems $\S=\{X_\alpha,p_\alpha^\beta,\A\}$ over a directed
set $\A$. Let us also note that for any continuous functor $F:Comp \to Comp$
the map $F:C(X,Y)\to C(FX,FY)$ (the space
$C(X,Y)$ is considered with the compact-open topology) is continuous.

Finally, a functor $F$ is called {\it weight-preserving} if $w(X)=w(F(X))$ for
every infinite $X\in Comp$.

A functor $F$ is called {\it normal} \cite{Shchepin} if it is
continuous, monomorphic, epimorphic, preserves
weight,intersections,preimages,singletons and the empty space. A
functor $F$ is said to be {\it weakly normal} ({\it almost
normal}) if it satisfies all the properties from the definition of
a normal functor excepting perhaps the preimage-preserving
property (epimorphness)(see \cite{ZT} for more details).

Similarly, one can define the same properties for a functor
$F:Tych\to Tych$ with the only difference that the
property of preserving surjections is replaced by the property of
sending $k$-covering maps to surjections (recall that $f:X\to Y$
is a $k$-covering map if for any compact set $B\subset Y$ there
exists a compact set $A\subset X$ with $f(A) = B$) (see \cite{ZT}, Def.2.7.1).

A.Chigogidze defined an extension construction of a functor in
$Comp$ onto $Tych$ the following way \cite{ChExt}. For any normal
functor $F:Comp\to Comp$ and any $X\in Tych$ the space $F_\beta(X)
= \{a\in F(\beta X)|\mathrm{there \ exists \ a \ compact \ set \
A\subset X \ with \ } a\in F(A)\}$ is considered with the topology
induced from $F(\beta X)$, where $\beta X$ is the Stone-Cech
compactification of the space $X$. Next, given any continuous
mapping $f:X\to Y$ between Tychonov spaces, put $F_\beta(f) =
F(\beta f)|_{F_\beta(X)}$. Then $F_\beta$ forms a covariant
functor in the category $Tych$. Chigogidze showed that in case $F$
is normal, the functor $F_\beta$ is also normal.

{\bf 2.} Let us modify the Chigogidze construction for any
functor $F:Comp\to Comp$. For $X\in Tych$ we put $F_\beta(X) =
\{a\in F(\beta X)|\mathrm{there \ exists \ a \ compact \ set \
A\subset X \ with \ } a\in F(i_A)(F(A))\}$ where by $i_A$ we
denote the natural embedding $i_A:A\hookrightarrow X$ (we do not
assume that the map $F(i_A)$  is an embedding). Evidently
$F_\beta$ preserves empty set and one-point space iff $F$ does.

Now we consider the problem when $F_\beta$ preserves embeddings.
Extension of any normal functor preserves embeddings, but, if we
drop the preimage preserving property, the situation could be
different. However, the examples from the introduction show that
the preimage-preserving property is not necessary. We define some
weaker property which will give us a necessary and sufficient
condition.

\begin{def1}We say that a monomorphic functor $F:Comp\to Comp$ \emph{preserves 1-preimages}, if
for any $f:X\to Y$, where $X,Y\in Comp$, any closed $A\subset Y$
such that $f|_{f^{-1}(A)}$ is a homeomorphism, we have that
$(Ff)^{-1}(FA) = F(f^{-1}(A))$. (Let us remark it is equivalent to
the condition that the map $Ff\mid(Ff)^{-1}(FA)$ is a
homeomorphism.)\end{def1}

Let us note that this definition was independently introduced by
T.Banakh and A.Kucharski \cite{BK}.

\begin{Pr1} If $F$ is a monomorphic functor that preserves
1-preimages in the class of open mappings, then $F$ preserves
1-preimages.\end{Pr1}

\begin{proof} Take any mapping $f:X\to Y$ such that
$f|_{f^{-1}(A)}$ is a homeomorphism for some closed subset
$A\subset Y$. Let $i_1:X\to X\times Y$ be an embedding defined by
the formula $i_1(x) = (x,f(x))$. Denote $Z = X\times Y /
\varepsilon$, where the relation $\varepsilon$ is given by
$\varepsilon = \{pr_Y^{-1}(a)|a\in A\}$ ($pr_Y:X\times Y$ is the
respective projection). Let $q:X\times Y\to Z$ be the quotient
mapping. The map $h:Z\to Y$ given by the conditions $h(z) = y$ for
any $z = (x,y)\in Z\setminus q(X\times A)$ and $h(z) = a$ for any
$z = q(pr_Y^{-1}(a))$, $a\in A$, is open and satisfies the
following two conditions: $pr_Y = h\circ q$, $h|_{h^{-1}(A)}$ is a
homeomorphism. Apparently, the map $i = q\circ i_1$ is an
embedding, moreover, $h\circ i = f$. Since $F$ preserves
1-preimages in the class of open mappings, we have $(Fh)^{-1}(FA)
= F(h^{-1}(A))$, which gives us the equality $(Ff)^{-1}(FA) =
F(f^{-1}(A))$.
\end{proof}

\begin{Pr2} If $F$ is a monomorphic functor that preserves
1-preimages, then $F_\beta$ preserves embeddings.\end{Pr2}

\begin{proof} Take any embedding $f:X\to Y$. Then the map $F_\beta(f)$ is closed as the restriction
of a closed map onto a full preimage and, moreover, injective,
hence an embedding.\end{proof}

For any $X\in Tych$ and any its compactification $bX$ we can
define $F_b(X) = \{a\in F(bX)|$ there is a compact subset
$A\subset X$ with $a\in F(A)\}\subset F(bX)$ and consider it with
the respective subspace topology.

\begin{Cor} If $F$ is a monomorphic, 1-preimage-preserving functor, then $F_\beta(X)\cong F_b(X)$ for
any Tychonov space $X$ and its compactification $bX$.\end{Cor}

\begin{Pr3} If $F$ is monomorphic, preserves 1-preimages and weight, then $F_\beta$ preserves
weight.\end{Pr3}

\begin{proof} The statement follows from the previous corollary
and the fact that for any $X\in Tych$ there exists its
compactification $bX$  which has the same weight as
$X$.\end{proof}

As the following proposition shows, the reverse implication to
that of Proposition 2 also holds.

\begin{Pr4} Let $F$ be a continuous functor such that $F_\beta$
preserves embeddings. Then $F$ preserves 1-preimages. \end{Pr4}

\begin{proof} Assume the contrary. Then there exist a map $f:X\to Y$ and a closed subset $A\subset Y$
such that $f|_{f^{-1}(A)}$ is a homeomorphism and $Ff^{-1}(FA)
\neq F(f^{-1}(A))$. Hence we can choose $\nu\in FA$ and $\mu\in
FX\backslash F(f^{-1}(A))$ such that $Ff(\mu) = \nu$. We will
construct a space $S\in Tych$ and its compactification $\gamma S$
such that the map $F_\beta(id_S):F_\beta(S)\to F_\beta(\gamma
S)=F(\gamma S)$ is not an embedding, where $id_S:S\to (\gamma S)$
is an identity embedding.

First put $Z = X\times \alpha \mathbb{N}$, where the space of
natural numbers $\mathbb{N}$ is considered with the discrete
topology and $\alpha \mathbb{N} = \mathbb{N}\cup \{\xi\}$ is the
one-point compactification of $\mathbb{N}$.  Define a continuous
function $g:Z\to Y$ by $g(x,n) = f(x)$ for any $x\in X, \ n\in
\alpha \mathbb{N}$. Let $T = Z/\varepsilon$ be a quotient space,
where $\varepsilon$ is an equivalence relation defined by its
classes of equivalence $\{\{x\}|x\in (X\setminus A)\times
\mathbb{N}\} \cup \{g^{-1}(y)\cap X\times\{\xi\}|y\in Y\setminus
A\}\cup\{\{a\} \times\alpha \mathbb{N}|a\in  A\} $. By $q:Z\to T$
we denote the respective quotient mapping. Then the map $h:T\to Y$
defined by the equality $g=h\circ q$ is continuous. The set $D =
q(X\times \{\xi\})$ is compact as a continuous image of a compact
set and moreover $h|_D$ is one-to-one, hence a homeomorphism
between $D$ and $Y$. We denote by $j:Y\to T$ the inverse
embedding. Also, for any $n\in \mathbb{N}$ the space $S_n =
q(X\times \{n\})$ is homeomorphic to $X$ and we denote by
$j_n:X\to T$ the inverse embedding. Then we have $h\circ j_n = f$.
Finally note that $T$ is a compactification of the space $S =
T\backslash  q((X\setminus A)\times\{\xi\})$.

Put $\mu_n=F(j_n)(\mu)$ for  $n\in \mathbb{N}$. The sequence $j_n$
converges to $j\circ f$ in the space $C(X,T)$. Since $F$ is
continuous, the sequence  $F(j_n)$ converges to $F(j\circ f)$ in
the space  $C(FX,FT)$. Hence the sequence $\mu_n$ converges to
$F(j\circ f)(\mu)=F(j)(\nu)\in F(q(A\times \alpha\mathbb{N}))$.

Now consider $F_\beta (S)$ as a subspace of $F(\beta S)$. There
exists a map $s_1:S\to X$ such that $s_1\circ j_n=\id_X$. Let
$s:\beta S\to X$ be the extension of $s_1$. Then
$Fs(\mu_n)=\mu\notin F(f^{-1}(A))$. Then the sequence $\mu_n$ does
not converge to any element of $F(q(A\times \alpha\mathbb{N}))$.
The proposition is proved.
\end{proof}

Propositions 2 and 4 yield the following

\begin{th1} For any continuous monomorphic functor $F$  the functor
$F_\beta$ preserves embeddings if and only if $F$ preserves
1-preimages.\end{th1}

The proof of the following proposition is a routine checking and we omit it.

\begin{Pr5}Let $F:Comp\to Comp$ be a functor.

1) if $F$ preserves embeddings, 1-preimages and intersections then
$F_\beta$ preserves intersections;

1) if $F$ preserves embeddings and preimages then
$F_\beta$ preserves preimages;

3) if $F$ preserves surjections then $F_\beta$ sends $k$-covering
maps to surjections;
\end{Pr5}

Now let us consider continuity of the Chigogidze extension. The
following example shows that in the absence of the
preimage-preserving property of the functor $F$, it is difficult
to speak of continuity of $F_\beta$, since even the extension of
such known weakly normal functor as $G$ does not possess it.

{\bf Example.} Let us define the inclusion hyperspace functor $G$.
Recall that a closed subset ${\cal A}\in \exp^2 X$, where $X\in Comp$ is
called an inclusion hyperspace, if for every $A\in{\cal A}$
and every $B\in \exp X$ the inclusion $A\subset B$ implies $B\in {\cal
A}$. Then $GX$ is the space of all inclusion hyperspaces
with the induced from $\exp^2 X$ topology. For any map $f:X\to Y$ define
$Gf:GX\to GY$ by $Gf({\cal A}) = \{B\in \exp Y| f(A) \subset B \ for \ some \ A\in{\cal A}\}$.
The functor $G$ is weakly normal (see \cite{ZT} for more details).
In the next section we will see that the functor $G$ preserves 1-preimages.

Let us show that the functor $G_\beta$ is not continuous.
Consider the following inverse system. For any $n\in
\mathbb{N}$ put $X_n = \mathbb{N}\times \{1,...,n\}$ (here the
spaces $\mathbb{N}$ and $\{1,...,n\}$ are considered with the
discrete topology). Define $p^m_n:X_m\to X_n$, where $m\ge n$, the
following way: $p^m_n(x,k) = (x,\min \{k,n\})$. We obtained the
inverse system $S = \{X_m,p^m_n,\mathbb{N}\}$. Then the limit
space $X = \lim S$ is homeomorphic to the space $\mathbb{N}\times
A$ (here $A = \alpha\mathbb{N} = \mathbb{N}\cup \{\xi\}$ is the
one-point compactification of $\mathbb{N}$, i.e. a convergent
sequence; also we put $\xi$ to be greater than any natural
number), and the limit projections $p_n:X\to X_n$ can be given by
$p_n(x,k) = (x,\min\{k,n\})$, $k\in A$. The continuity of
$G_\beta$ means that $\lim G_\beta(p_n):G_\beta(\lim S)\to \lim
G_\beta(S)$ is a homeomorphism. Here both $G_\beta(\lim S)$ and
$\lim G_\beta(S)$ can be thought as subspaces of $G(bX)$, where
$b$ is a compactification of $X$ with the property $bX = \lim
\beta S$.

Now we will construct $K\in \lim G_\beta(S)$ which does not belong
to $\lim G_\beta(p_n)(G_\beta(\lim S))$. Consider the space $X$
imbedded into its compactification $bX$. For any $n\in A\backslash
\{\xi\}$ put $K_n = \{1,...,n\}\times\{n\}$. If we want to obtain
a closed family of sets, the set $K_\xi =
\overline{\mathbb{N}\times\{\xi\}}$ must be added to the family
$\widetilde{K} = \{K_n\}_{n\in \mathbb{N}}$. Now put $K =
\{B\subset bX|K_n\subset B$ for some $n\in A\}$. Then $K\in \lim
G_\beta(S)$. However, there is apparently no element $C\in
G_\beta(\lim S)$ with $\lim G_\beta(p_n)(C) = K$. Hence, $\lim
G_\beta(p_n)$, being not surjective, is not a homeomorphism.

{\bf 3.} We start this section  with definitions of some functors
we deal with in this paper. Let $X$ be compactum. By $C(X)$ we
denote the Banach space of all continuous functions $\phi:X\to\R$
with the usual $\sup$-norm. We consider $C(X)$ with natural order.
Let $\nu :C(X)\to\mathbb{R}$ be a functional (we do not suppose a
priori that $\nu$ is linear or continuous). We say that $\nu$ is
1) non-expanding if $\vert \nu(\varphi)-\nu(\psi)\vert\le
d(\varphi,\psi)$ for all $\varphi,\psi\in C(X)$;  2) weakly
additive if for any function $\phi \in C(X)$ and any $c\in
\mathbb{R}$ we have $\nu(\phi+c_X)=\nu(\phi)+c$ (by $c_X$ we
denote the constant function); 3) preserves order if for any
$\varphi,\psi\in C(X)$ such that $\varphi\le \psi$ the inequality
$\nu(\varphi)\le\nu(\psi)$ holds; 4) linear if for any $\alpha$,
$\beta\in\mathbb{R}$ and for any two functions $\psi$, $\phi \in
C(X)$ we have
$\nu(\alpha\phi+\beta\psi)=\alpha\nu(\phi)+\beta\nu(\psi)$.

Now for any space $X$ denote $VX = \prod_{\varphi\in
C(X)}[\min\varphi,\max\varphi]$. For any mapping $f:X\to Y$ define
the map $Vf$ as follows: $Vf(\nu)(\varphi) = \nu(\varphi\circ
f)$ for every $\nu\in VX, \ \varphi\in C(Y)$. Then $V$ is a
covariant functor in the category $Comp$ \cite{OnV}.

Let us remark that the space $VX$ could be considered as the space
of all functionals $\nu :C(X)\to\mathbb{R}$ with the only
condition $\min\varphi(X)\le\nu(\varphi)\le \max\varphi(X)$ for
every $\nu\in VX, \ \varphi\in C(Y)$. By $EX$ we denote the subset
of $VX$ defined by the condition 1) (non-expanding functionals;
see \cite{Ca} for more details), by $EAX$ the subset defined by
the conditions 1) and 2). The conditions 2) and 3) define the
subset $OX$ (order-preserving functionals, see \cite{RO});
finally, the conditions 3) and 4) define the well-known subset
$PX$ (probability measures, see for example  \cite{PF}). For a map
$f:X\to Y$ the mapping $Ff$, where $F$ is one of $P$, $O$, $EA$,
$E$, is defined as the restriction of $Vf$ on $FX$. It is easy to
check that the constructions $P$, $O$, $EA$ and $E$ define
subfunctors of $V$. It is known that the functors $O$ and $E$ are
weakly normal (see \cite{RO} and \cite{Ca}). Using the same
arguments one can check that  $EA$ is weakly normal too.

The question arises naturally which of defined above functors have
the property of preserving 1-preimages. It is easy to check that
we have the inclusions $PX\subset OX\subset EAX\subset EX\subset
VX$. We will show that the functor $EA$ satisfies this property
and $E$ does not. Since subfunctors inherit the 1-preimages
preserving property, this is the complete answer. Let us also
remark that the results of  \cite{OnV} and \cite{LM} show that
many other known functors could be considered as subfunctors of
$EA$, for example the superextension, the hyperspace functor, the
inclusion hyperspace functor etc. This shows that the class of
functors with the 1-preimages preserving property is wide enough.

We start with a definition of an $AR$-compactum. Recall that a
compactum $X$ is called an absolute retract (briefly $X\in AR$) if
for any embedding $i:X\to Z$ of $X$ into compactum $Z$ the image
$i(X)$ is a retract of $Z$.

The next lemma will be needed in the following discussion.

\begin{Le1} Let $F$ be a monomorphic subfunctor of $V$ which preserves
intersections and $B$ be a closed subset of a compactum $X$. Then
$\nu\in FB$ iff $\nu(\varphi_1)=\nu(\varphi_2)$ for each
$\varphi_1$, $\varphi_2\in C(X)$ such that $\varphi_1|_B=\varphi_2|_B$.\end{Le1}

\begin{proof} \emph{Necessity}. The inclusion $\nu\in FB\subset FX$ means
that there exists $\nu_0\in FB$ with $F(i_B)(\nu_0) = \nu$, where
$i_B:B\to X$ is a natural embedding. Hence, for any $\varphi_1$,
$\varphi_2\in C(X)$ such that $\varphi_1|_B=\varphi_2|_B$ we have
$\nu(\varphi_1) = \nu_0(\varphi_1\circ i_B) = \nu_0(\varphi_2\circ
i_B) = \nu(\varphi_2)$.

\emph{Sufficiency}. We can find an embedding $j:B\hookrightarrow
Y$, where $Y\in AR$. Define $Z$ to be the quotient space of the
disjoint union $X\cup Y$ obtained by attaching $X$ and $Y$ by $B$.
Denote by $r:Z\to Y$ the retraction mapping.

Now take any $\nu\in FX\subset FZ$ with the property
$\nu(\varphi_1)=\nu(\varphi_2)$ for each $\varphi_1$,
$\varphi_2\in C(X)$ such that $\varphi_1|_B=\varphi_2|_B$. We
claim that $F(r)(\nu) = \nu$. Indeed, take any $\varphi\in C(Z)$.
Then $F(r)(\nu)(\varphi) = \nu(\varphi\circ r) = \nu(\varphi)$
since $\varphi\circ r|_Y = \varphi|_Y$. Hence, $\nu\in FX\cap FY =
FB$.
\end{proof}

\begin{Pr6} The functor $EA$ preserves 1-preimages.\end{Pr6}

\begin{proof} Let $f:X\to Y$ be a continuous open map between compacta $X$ and
$Y$ and $B$ be a closed subset of $Y$ such that $f|_{f^{-1}(B)}$
is a homeomorphism. Choose any $\nu\in EA(B)\subset EA(Y)$. Using
Lemma 1 we can define $\mu_0\in EA(f^{-1}(B))$ by the condition
$\mu_0(\varphi)=\nu(\psi)$ for each $\varphi\in C(X)$ and $\psi\in
C(Y)$ such that $\psi\circ f|f^{-1}(B)=\varphi|_{f^{-1}(B)}$.

It is enough to show that for each $\mu\in(EA(f))^{-1}(\nu)$ we
have $\mu=\mu_0$. Suppose the contrary. Then there exist
$\varphi\in C(X)$ and $\psi\in C(Y)$ such that $\psi\circ
f|f^{-1}(B)=\varphi|_{f^{-1}(B)}$ and $\mu(\varphi)\ne\nu(\psi)$.
We can suppose that $\mu(\varphi)>\nu(\psi)$. Define a function
$\psi':Y\to \mathbb{R}$ by $\psi'(y) = \max \varphi f^{-1}(y)$ for
any $y\in Y$. The function $\psi'$ is continuous since $f$ is
open. Put $\xi=(\psi' - D)\circ f$, where $D = \sup\{\max \varphi
f^{-1}(y)- \min \varphi f^{-1}(y)|y\in Y\}$. Then
$d(\xi,\varphi)\le D$ but
$\mu(\varphi)-\mu(\xi)=\mu(\varphi)-\mu((\psi' - D)\circ
f)=\mu(\varphi)-\nu(\psi') + D=\mu(\varphi)-\nu(\psi) + D>D$ and
we obtain a contradiction. The proof is similar for the case
$\mu(\varphi)<\nu(\psi)$.

Hence, $EA$ preserves 1-preimages in the class of open mappings,
and, by Proposition 1, we are done.
\end{proof}

\begin{Pr7} The functor of nonexpanding functionals $E$ does not
preserve 1-preimages.\end{Pr7}

\begin{proof} Consider the mapping $f:X\to Y$ between
discrete spaces $X =\{x,y,s,t\}$ and $Y = \{a,b,c\}$ which is
defined as follows: $f(x) = a$, $f(y) = b$, $f(s) = f(t) = c$. Put
$A = \{\varphi\in C(X)|\varphi(s) = \varphi(t)\}$. Define the
functional $\nu:A\to \mathbb{R}$ as follows: $\nu(\varphi) = \min
\{\varphi(x),\varphi(y)\}$ if $\varphi|_{\{x,y\}}\ge 0$,
$\nu(\varphi) = \max \{\varphi(x),\varphi(y)\}$ if
$\varphi|_{\{x,y\}}\le 0$, and $\nu(\varphi) = 0$ otherwise. One
can check that $\nu$ is nonexpanding. Now take the function
$\psi:X\to \mathbb{R}$ defined as follows $\psi(x) =
1$, $\psi(y) = -1$, $\psi(s) = 0$, $\psi(t) = 4$.
One can check that we can extend $\nu$ to a nonexpanding
functional on $A\cup \{\psi\}$ by defining its value on $\psi$ to
be  $-1$. This new functional can be further extended
to a nonexpanding functional on the whole $C(X)$ \cite{Ca}. Denote this
extension by $\widetilde{\nu}$. Evidently, $Ef(\widetilde{\nu})\in
E(\{a,b\})$. On the other hand, $\widetilde{\nu}\notin
E(\{x,y\})$.
\end{proof}

{\bf 4.} We consider in this section a monomorphic continuous
functor $F$ which preserves intersections, weight, empty set,
point and 1-preimages. We investigate topology of the space
$F_\beta Y$ where $Y$ is a metrizable separable non-compact space.
We consider $Y$ as a dense subset of metrizable compactum $X$. It
follows from Corollary 1 that $F_\beta Y$ is homeomorphic to
$F_bY\subset FX$ (where $X$ is considered as a compactification
$bY$ of $Y$) and in what follows we identify $F_\beta Y$ with
$F_bY$. Also, the properties we impose on $F$ imply that $F_\beta
Y$ is a dense proper subspace of $FX$.

T.Banakh proved in \cite{BB} that  $F_\beta Y$ is
$F_\sigma$-subset of $FX$ when $Y$ is locally compact; $F_\beta Y$
is $F_{\sigma\delta}$-subset when $Y$ is $G_\delta$-subset. If $Y$
is not a $G_\delta$-subset, then $F_\beta Y$ is not analytic.

We consider in the Hilbert cube $Q = [-1,1]^\omega$ the following
subsets: $\Sigma = \{(t_i)\in Q|\sup_i|t_i|<1\}$; $\sigma =
\{(t_i)\in Q| t_i\neq 0 \mathrm{for \ finitely \ many \ of } \
i\}$ and $\Sigma^\omega\subset Q^\omega\cong Q$.

It is shown in \cite{PBC} that any analytic $P_\beta Y$ is
homeomorphic to one of the spaces $\sigma$, $\Sigma$ or
$\Sigma^\omega$. We generalize this result for convex functors.

By $Conv$ we denote the category of convex compacta (compact
convex subsets of locally convex topological linear spaces)  and
affine maps. Let $U:Conv\to Comp$ be the forgetful functor. A
functor $F$ is called {\it convex} if there exists a functor
$F':Comp\to Conv$ such that $F=UF'$. It is easy to see that the
functors $V$, $E$, $EA$, $O$ and $P$ are convex. It is shown in
\cite{Sha} that for each convex functor $F$ there exists a unique
natural transformation $l:P\to F$ such that the map $lX:PX\to FX$
is an affine embedding.

\begin{Le2} $P_\beta Y=(lX)^{-1}(F_\beta Y)$. \end{Le2}

\begin{proof} Take any measure $\mu\in P(X)$ such that $lX(\mu) = \mu'\in F_\beta
Y$. By the definition of $F_\beta Y$ it means that $\mu'\in FB$
for some compactum $B\subset Y$. We will show that $\mu\in PB\ss
P_\beta Y$. Choose an absolute retract $T$ which contains $B$ and
define $Z$ to be the quotient space of the disjoint union $X\cup
T$ obtained by attaching $X$ and $T$ by $B$. By $r:Z\to T$ denote
the retraction. Since $l$ is a natural transformation and $r$ is
an identity on $T\subset Z$, we have that $F(r)\circ lZ(\mu) =
\mu' = lT\circ P(r) (\mu)$. Hence, $\mu = P(r)(\mu)\in P(T)$ due
to injectivity of $lZ$. Therefore, $\mu\in PX\cap PT = PB$. The
lemma is proved.
\end{proof}

We need some notions from infinite-dimensional topology. See
\cite{BRZ} for more details. All spaces are assumed to be
metrizable and separable. A closed subset $A$ of a compactum $T$
is called $Z$-set if there exists a homotopy $H:T\times [0;1]\to
T$ such that $H\mid_{T\times\{0\}}=\id_{T\times\{0\}}$ and
$H(T\times (0,1])\cap A=\emptyset$; a subset $B$ of $T$ is called
$\sigma Z$-set if it is contained in countable union of $Z$-sets
of $T$. In what follows we will use the following facts.

We don't know if $F_\beta Y$ is a $\sigma Z$-set in $FX$ for any
convex functor $F$. Thus, we introduce some additional property.
We consider the compactum $FX$ as a convex subset of a locally
convex linear space.

\begin{def2} A convex functor $F:Comp\to Comp$ is called strongly
convex if for each compactum $X$, each closed subset $A\subset X$
we have $(FX\setminus FA)\cap \mathrm{aff} FA =
\emptyset$.\end{def2}

\begin{Pr8} Each convex subfunctor $F$ of the functor $V$ is strongly convex.\end{Pr8}

\begin{proof} By Lemma 1 any element from $\mathrm{aff} FA$ takes
the same value at any two functions from $C(X)$ which coincide on
$A$, which is not true for functionals from $FX\setminus FA$.
\end{proof}

\begin{Pr9} Let $F$ be a strongly convex functor. Then $F_\beta Y$ is a $\sigma Z$-set in $FX$.\end{Pr9}

\begin{proof}  Take any $y\in X\backslash Y$. Then $F_\beta Y\subset
F_\beta(X\backslash \{y\})$, and $X\backslash\{y\}$ can be
represented as a countable union of its compact subsets $A_n$ with
the property that $A_n\subset \mathrm{int} A_{n+1}$, hence,
$F_\beta(X\backslash\{y\}) = \cup_{n\in\mathbb{N}}F(A_n)$. Let us
show that all $F(A_n)$ are $Z$-sets in $FX$. Take any $\nu\in
FX\setminus F_\beta(X\setminus\{y\})$ and the set $Z =
\{t\nu+(1-t)\mu|t\in(0,1], \mu\in F_\beta (X\setminus\{y\})\}$.
Since $F$ is strongly convex, we have $Z\cap F_\beta
(X\setminus\{y\}) = \emptyset$. Since $Z$ is convex and dense
subset of $FX$, there exists a homotopy $H:FX\times[0,1]\to FX$
such that $H(FX\times(0,1])\subset Z$ (see, for example, Ex. 12,
13 to section 1.2 in \cite{BRZ}).
\end{proof}

Now, we are going to obtain the complete topological
classification of the pair $(FX,F_\beta Y)$ where $X$ is a
metrizable compactum and $Y$ its proper dense $G_\delta$-subset.
We need some characterization theorems.

\begin{thA}\cite{D} Let $C$ be an infinite-dimensional dense $\sigma Z$ convex subspace of a
a convex metrizable compactum $K$, and additionally let $C$ be a
countable union of its finite-dimensional compact subspaces. Then
the pair $(K,C)$ is homeomorphic to $(Q,\sigma)$.
\end{thA}

\begin{thB}\cite{CDM} Let $K$ be a convex metrizable compactum, and let $C\subset K$ be
its proper dense $\sigma Z$ convex $\sigma$-compact subspace that
contains an infinite-dimensional convex compactum. Then the pair
$(K,C)$ is homeomorphic to the pair $(Q,\Sigma)$.
\end{thB}

The following theorem follows from 5.3.6, 5.2.6, 3.1.10 in
\cite{BRZ}.

\begin{thC} Let $K$ be a convex compact subset locally convex linear metric space, and let $C\subset K$ be
its proper dense $\sigma Z$ convex $F_{\sigma\delta}$ subspace
such that $K\setminus C)\cap \mathrm{aff} C = \emptyset$, and
additionally there exists a continuous embedding $h:Q\to K$ such
that $h^{-1}(C)=\Sigma^\omega$. Then the pair $(K,C)$ is
homeomorphic to the pair $(Q,\Sigma^\omega)$.
\end{thC}

\begin{TopPr} Let $F$ be a strongly convex functor, $X$ is a metrizable compactum and $Y$ is its proper dense $G_\delta$-subset.
The pair $(FX,F_\beta Y)$ is homeomorphic to

\begin{enumerate}

\item $(Q,\sigma)$, if $Y$ is discrete subspace of $X$ and $F(n)$
is finite-dimensional for each $n\in\N$;

\item $(Q,\Sigma)$, if $Y$ is discrete subspace of $X$ and $F(n)$
is infinite-dimensional for some $n\in\N$ or $Y$ is locally
compact non-discrete subspace of $X$;

\item $(Q,\Sigma^\omega)$, if $Y$ is not locally compact.
\end{enumerate} \end{TopPr}

\begin{proof} It is easy to see that $F_\beta Y$ is a convex subset of $FX$.

We prove the first assertion. Since $X$ is metrizable, $Y$ is
countable. We can represent $Y=\cup_{n=1}^\infty Y_n$ where
$|Y_n|=n$. Then $F_\beta Y= \cup_{n=1}^\infty FY_n$. Since $PY_n$
could be considered as an $n-1$-dimensional subspace of $FY_n$,
the space $F_\beta Y$ is infinite-dimensional. Moreover,  $F_\beta
Y$ is a $\sigma Z$-set by Proposition 9. Since each $FY_n$ is a
finite-dimensional compactum, we can apply Theorem A.

We prove the second assertion. In the case when $Y$ is discrete,
$FY_n$ is infinite-dimensional convex compactum for some $n$. When
$Y$ is not discrete, it contains an infinite compactum $Y'$ and
$FY'$ is infinite-dimensional convex compactum. We apply
Proposition 9 and Theorem B.

For the third assertion, note that the pair $(PX,P_\beta Y)$ is
homeomorphic to $(Q,\Sigma^\omega)$ \cite{PBC}. Since $F$ is
strongly convex, we have $(FX\setminus F_\beta Y)\cap \mathrm{aff}
F_\beta Y = \emptyset$. We apply Lemma 2, Proposition 9 and
Theorem C.
\end{proof}

\begin{Cor2} Suppose that $F$ is a strongly convex functor.
Then for any separable metrizable space $X$

1) $X\cong \mathbb{N}$ implies $F_\beta(X)\cong Q_f$ in case
$F(n)$ is finite-dimensional for any $n\in\mathbb{N}$ or
$F_\beta(X)\cong \Sigma$ otherwise;

2) if $X$ is locally compact non-discrete and non-compact then
$F_\beta(X)\cong \Sigma$;

3) if $X$ is topologically complete not locally compact then
$F_\beta(X)\cong \Sigma^\omega$. \end{Cor2}


\begin{thebibliography}{100}

\bibitem{BB} T.~Banakh, {\it Descriptive classes of sets and topological functors}, Ukrain. Mat. Zh., {\bf 47} (1995) 408--410.

\bibitem{PBC} T.~Banakh, R.~Cauty, {\it Topological classification
of spaces of probability measures over coanalytic sets},  Mat.
Zametki. \textbf{55} (1994), 10–19 (Russian).

\bibitem{BK} T.~Banakh, M.~Klymenko, A.~Kucharski, {\it On functors preserving skeletal maps and skeletally generated compacta}, (in preparation).

\bibitem{BRZ} T.~Banakh, T.~Radul, M.~Zarichnyi, Absorbing sets in
infinite-dimensional manifolds, VNTL Publishers. Lviv, 1996.

\bibitem{Ca} J.~Camargo, {\it The functor of nonexpanding functionals.}, Rev. Integr. Temas Mat., {\bf 20} (2002) 1--12.

\bibitem{ChExt} A.~Chigogidze, {\it On extension of normal
functors}, Vestnik Mosk. univ. Mat. Mekh. {\bf 6}, 23-26
(Russian).

\bibitem{CDM} D.W.~Curtis, T.~Dobrowolsky, J.~Mogilski, {\it Some
applications of the topological characterizations of the
sigma-compact spaces $l_2^f$ and $\Sigma$}, Trans. Amer. Math.
Soc., {\bf 284}, 837--846.

\bibitem{D} T.Dobrowolsky {\it The compact Z-property in convex
sets}, Top.Appl., {\bf 23}, 163--172.

\bibitem{FedZar} V.~Fedorchuk, M.~Zarichnyi.,  Covariant functors in
categories of topological spaces,  Results of Science and
Thechnics. Algebra.Topology.Geometry. Moscow.VINITI, v.28
P.47--95.

\bibitem{RO} T.~Radul, {\it On the functor of order-preserving functionals}, Commentat. Math. Univ.  Carol., {\bf 39} (1998) 609--615.

\bibitem{OnV} T.~Radul, {\it On strongly Lawson and I-Lawson monads}, Boletin de Mathematicas, {\bf 6} (1999), 69--76.

\bibitem{LM} T.~Radul, {\it On functional representations of
Lawson monads}, Applied Categorical Structures {\bf 9} (2001)
69--76.

\bibitem{ZT} A.~Teleiko, M.~Zarichnyi, Categorical Topology of
Compact Hausdorff Spaces, VNTL Publishers. Lviv, 1999.

\bibitem{Sha} L.~Shapiro, {\it On operators of extension of functions and normal functors}, Vestn. Mosk. univ., {\bf 1} (1992) 35--42.

\bibitem{Shchepin} E.~Schepin, {\it Functors and uncountable powers of compacta}, Uspekhi Mat. Nauk, {\bf 36} (1981),
3--62 (Russian).

\end{thebibliography}
\end{document}